\documentclass[11pt]{article}

\usepackage{amsmath, amssymb, amscd}

\def\eqref#1{(\ref{#1})}
\newcommand{\goth}{\frak}
\newcommand{\g}{{\frak g}}
\newcommand{\arrow}{{\:\longrightarrow\:}}

\newcommand{\C}{{\Bbb C}}
\newcommand{\R}{{\Bbb R}}

\renewcommand{\H}{{\Bbb H}}
\newcommand{\6}{\partial}
\def\1{\sqrt{-1}\:}

\renewcommand{\bar}{\overline}
\renewcommand{\phi}{\varphi}
\renewcommand{\epsilon}{\varepsilon}
\renewcommand{\geq}{\geqslant}
\renewcommand{\leq}{\leqslant}

\newcommand{\End}{\operatorname{End}}

\newcommand{\Hol}{\operatorname{Hol}}

\newcommand{\comment}[1]{{}}

\def\blacksquare{\hbox{\vrule width 4pt height 4pt depth 0pt}}
\def\endproof{\blacksquare}

\makeatletter

\@ifundefined{Bbb}
     {\newcommand{\Bbb}[1]{{\mathbb #1}}}%
{}

\newcommand{\ps@verbit}{%
  \renewcommand{\@oddhead}{%
          \scriptsize
          {Hypercomplex manifolds with trivial canonical bundle}
          \hfil\tiny {M. Verbitsky, June 23, 2004}}
  \renewcommand{\@evenhead}{\@oddhead}
  \renewcommand{\@oddfoot}{\hfil\thepage\hfil}
  \renewcommand{\@evenfoot}{\@oddfoot}}
 
\newcounter{Mycounter}[section]
\newcounter{lemma}[section]
\renewcommand{\thelemma}{{Lemma \thesection.\arabic{lemma}}}
\newcommand{\lemma}{%
     \setcounter{lemma}{\value{Mycounter}}
     \refstepcounter{lemma}
     \stepcounter{Mycounter}
     {\bf \thelemma:\ }}

\newcounter{claim}[section]
\renewcommand{\theclaim}{{Claim \thesection.\arabic{claim}}}
\newcommand{\claim}{%
     \setcounter{claim}{\value{Mycounter}}
     \refstepcounter{claim}
     \stepcounter{Mycounter}
     {\bf \theclaim:\ }}

\newcounter{sublemma}[section]
\newcounter{corollary}[section]
\renewcommand{\thecorollary}{{Corollary \thesection.\arabic{corollary}}}
\newcommand{\corollary}{%
     \setcounter{corollary}{\value{Mycounter}}
     \refstepcounter{corollary}
     \stepcounter{Mycounter}
     {\bf \thecorollary:\ }}

\newcounter{theorem}[section]
\renewcommand{\thetheorem}{{Theorem \thesection.\arabic{theorem}}}
\newcommand{\theorem}{%
     \setcounter{theorem}{\value{Mycounter}}
     \refstepcounter{theorem}
     \stepcounter{Mycounter}
     {\bf \thetheorem:\ }}

\newcounter{conjecture}[section]
\newcounter{proposition}[section]
\renewcommand{\theproposition}
       {{Proposition \thesection.\arabic{proposition}}}
\newcommand{\proposition}{%
     \setcounter{proposition}{\value{Mycounter}}
     \refstepcounter{proposition}
     \stepcounter{Mycounter}
     {\bf \theproposition:\ }}

\newcounter{definition}[section]
\renewcommand{\thedefinition}
       {{Definition~\thesection.\arabic{definition}}}
\newcommand{\definition}{%
     \setcounter{definition}{\value{Mycounter}}
     \refstepcounter{definition}
     \stepcounter{Mycounter}
     {\bf \thedefinition:\ }}

\newcounter{example}[section]
\newcounter{remark}[section]
\renewcommand{\theremark}{{Remark \thesection.\arabic{remark}}}
\newcommand{\remark}{%
     \setcounter{remark}{\value{Mycounter}}
     \refstepcounter{remark}
     \stepcounter{Mycounter}
     {\bf \theremark:\ }}

\newcounter{problem}[section]
\newcounter{question}[section]
\renewcommand{\thequestion}{{Question \thesection.\arabic{question}}}
\newcommand{\question}{%
     \setcounter{question}{\value{Mycounter}}
     \refstepcounter{question}
     \stepcounter{Mycounter}
     {\bf \thequestion:\ }}

\@addtoreset{equation}{section}
\@addtoreset{footnote}{section}
\makeatother

\begin{document}

\begin{center}
{\LARGE\bf
Hypercomplex manifolds 
 with trivial  \\[3mm] canonical  bundle and their holonomy
}
\\[4mm]
Misha Verbitsky\footnote{ The author is
supported by EPSRC grant  GR/R77773/01 
and CRDF grant RM1-2354-MO02}
\\[4mm]

{\tt verbit@maths.gla.ac.uk, \ \  verbit@mccme.ru}
\end{center}

{\small 
\hspace{0.15\linewidth}
\begin{minipage}[t]{0.7\linewidth}
{\bf Abstract} \\
Let $(M,I,J,K)$ be a compact hypercomplex manifold
admitting a special kind of quaternionic-Hermitian 
metric called an HKT-metric. Assume that the canonical bundle
of $(M,I)$ is trivial as a holomorphic line bundle. We show
that the holonomy of canonical torsion-free connection called
Obata connection on $M$ is contained in $SL(n, {\Bbb H})$. In Appendix we
apply these arguments to compact nilmanifolds which admit
an abelian hypercomplex structuree, showing
that such manifolds have holonomy in $SL(n, {\Bbb H})$.
\end{minipage}
}

{
\small
\tableofcontents
}

\section{Hypercomplex manifolds and holonomy}

Let $(M, I, J, K)$ be a manifold equipped with an action of 
the quaternion algebra ${\Bbb H}$ on its tangent bundle. 
The manifold $M$ is called {\bf
hypercomplex} if the operators $I, J, K\in \Bbb H$ define
integrable complex structures on $M$. As Obata proved
(\cite{_Obata_}), this condition
is satisfied if and only if $M$ admits a torsion-free
connection $\nabla$ preserving the quaternionic action:
\[
\nabla I = \nabla J = \nabla K =0.
\]
Such a connection is called {\bf an Obata 
connection on $(M, I, J, K)$}. It is necessarily
unique (\cite{_Obata_}).

Hypercomplex manifolds were defined by C.P.Boyer
(\cite{_Boyer_}), who gave a classification of
compact hypercomplex manifolds for $\dim_{\Bbb H} M =1$.
Many interesting examples of hypercomplex manifolds
were found in 1990-ies, 
see e.g. \cite{_Joyce_}, \cite{_Poon_Pedersen_},
\cite{_Barberis_Dotti_}. An excellent 
survey and bibliography of this stage of 
research can be found in \cite{_Gra_Poon_}.

We are interested in the algebraic geometry of the complex
manifold $(M,I)$ underlying the hypercomplex structure.
In this paper we relate the holomorphic geometry of
the canonical bundle of $(M,I)$ with the differential
geometry of Obata connection. 

Let $\Hol(\nabla)$ be the holonomy group 
associated with the Obata connection $\nabla$.
Since $\nabla$ preserves the quaternionic structure,
$\Hol(\nabla)\subset GL(n, {\Bbb H})$. We define the
determinant of $h\in GL(n, {\Bbb H})$ in the following
way. Let $V\cong {\Bbb H}^n$ be the vector space over $\Bbb H$,
and $V_I^{1,0}$ the same space considered as a complex
space with the complex structure induced by $I$. 
The Hodge decomposition gives 
$V \otimes \C \cong V_I^{1,0}\oplus V_I^{0,1}$
The top exterrior power $\Lambda^{2n,0}_I(V):=\Lambda^{2n}(V_I^{1,0})\cong \C$
is equipped with a natural real structure:
\begin{equation}\label{_real_stru_from_J_Equation_}
\eta \arrow J(\bar\eta)
\end{equation}
for $\eta \in \Lambda^{2n,0}_I(V)$
(the quaternions $I$ and $J$ anticommute,
hence $J$ exchanges $\Lambda^{p,q}_I(V)$ with 
$\Lambda^{q,p}_I(V)$). Since the real structure
on $\Lambda^{2n,0}_I(V)$ is constructed from
the quaternion action, any $h\in GL(V,{\Bbb H})$
preserves this real structure. Let $\det(h)$
denote the action induced by $h$ on $\Lambda^{2n,0}_I(V)\cong \C$.
Then $\det (h)\in \R$, as the above argument imples.
This defines a homomorphism 
\[
\det:\; GL(n, {\Bbb H}) \arrow \R^*
\]
to the multiplicative group of non-zero real
numbers; clearly, $\det(h)$ is always positive.
Let $SL(n, \H)\subset GL(n, {\Bbb H})$ be the kernel of $\det$.

\hfill

From the above argument the following claim follows
immediately.

\hfill

\claim \label{_SL_holo_Claim_}
Let $(M, I, J, K)$ be a hypercomplex manifold,
$\dim_{\Bbb H}(M) =n$, $\nabla$ the Obata connection, 
$K(M,I)$ the canonical bundle, and $\nabla_K$ the connection
on $K(M,I)$ induced by the Obata connection. Then the
holonomy of $\nabla_K$ is trivial if and only if
$\Hol(\nabla)\subset SL(n, {\Bbb H})$.

\endproof

\hfill

The Obata connection $\nabla$ is torsion-free and preserves 
\[ \Lambda^{2n,0}_I(M)\subset \Lambda^{2n}(M).\] 
Therefore, the $(0, 1)$-part of $\nabla$
on $(M,I)$ is equal to the holomorphic structure operator on 
$(M,I)$.\footnote{The ``holomorphic structure operator'' on 
a vector bundle $B$ is understood as an operator $B\arrow \Lambda^{0,1}(M)$
vanishing on holomorphic sections and satisfying the Leibniz rule.}
This means that the $(0,1)$-part of $\nabla_K$
is the holomorphic structure operator on the
canonical bundle $K(M,I)$. This gives the following claim.

\hfill

\claim\label{_SL_n_H_tri_can_class_Claim_}
Let $(M, I, J, K)$ be a hypercomplex manifold,
and $\nabla$ its Obata connection.
Assume that \[ \Hol(\nabla)\subset SL(n, {\Bbb H}). \]
Then the canonical
bundle of $(M,I)$ is trivial.

\hfill

{\bf Proof:} By \ref{_SL_holo_Claim_}, the holonomy
of the Obata connection $\nabla_K$ on $K(M,I)$ is trivial. 
Let $\eta$ be a non-zero section of  $K(M,I)$
preserved by $\nabla_K$. Then $\nabla_K^{0,1}\eta=\bar\6\eta=0$,
hence $\eta$ is a holomorphic trivialization of $K(M,I)$.
\endproof 

\hfill

One can ask the following question.

\hfill

\question
Let $(M,I,J,K)$ be a compact hypercomplex manifold. 
Assume the complex manifold $(M,I)$ has trivial canonical
bundle. Does it follow that $\Hol(M)\subset SL(n,{\Bbb H})$?

\hfill

In this paper we give an affirmative answer to this question
(\ref{_K_triv=>SL(n,H)_Theorem_}),
provided that $M$ admits a special kind of quaternionic 
Hermitian metric,
\footnote{
A metric $g$ is called {\bf quaternionic Hermitian} if 
\[ g(Ix, Iy) = g(Jx, Jy) = g(Kx, Ky) = g(x, y)\] for all $x, y \in TM$.}
 so-called HKT-metric, or hyperk\"ahler
metric with torsion.

\hfill

In the Appendix, we show that any compact nilmanifold
admitting an abelian hypercomplex structure (\ref{_abelian_hc_Definition_})
has holonomy in $SL(n,{\Bbb H})$ 
(\ref{_SL(N,H)_nilma_Corollary_}). 
This gives plenty of examples of
hypercomplex non-\-hy\-per\-k\"ah\-ler
manifolds with (global) holonomy in 
$SL(n,{\Bbb H})$, though the local
holonomy in these cases is trivial.

\hfill

\remark 
It is easy to see that the canonical bundle is topologically
trivial for any hypercomplex manifold (see e.g. \cite{_Verbitsky:HC_Kahler_}).
On a compact K\"ahler manifold, topological triviality
of canonical bundle implies that it is trivial
holomorphically. This follows easily from the
Calabi-Yau theorem. On a non-K\"ahler manifold,
this is no longer true. In fact, for a typical non-hyperk\"ahler
compact hypercomplex manifold $(M, I, J, K)$,
the complex manifold $(M, I)$ admits no K\"ahler metrics,
and $K(M,I)$ is in most cases non-trivial
as a holomorphic vector bundle, though it is trivial
topologically. It is 
possible to show that $K(M,I)$ is non-trivial
for all hypercomplex manifolds $(M, I, J, K)$
such that $(M,I)$ is a principal toric fibration over
a base which has non-trivial canonical class; 
these include locally conformally hyperk\"ahler 
manifolds (see \cite{_Ornea:LCHK_})
and compact Lie groups with the hypercomplex structure
constructed by D. Joyce (\cite{_Joyce_}).

\section{HKT metrics on hypercomplex manifolds}
\label{_HKT_Subsection_}

Let $M$ be a hypercomplex manifold. A ``hyperk\"ahler with
torsion'' (HKT) metric on $M$ is a special
kind of a quaternionic Hermitian metric,
which became increasingly important in
mathematics and physics during the 
last 7 years.

HKT-metrics were
introduced by P.S.Howe and G.Papadopoulos (\cite{_Howe_Papado_})
and much discussed in the physics literature since then.
For an excellent survey of these works written from  a mathematician's
point of view, the reader is referred to the paper of
G. Grantcharov and Y. S. Poon \cite{_Gra_Poon_}.

The term ``hyperk\"ahler metric with torsion'' is actually
misleading, because an HKT-metric is not hyperk\"ahler.
This is why we prefer to use an abbreviation ``HKT-manifold''.

Let $(M, I, J, K)$ be a hypercomplex manifold,
$g$ a quaternionic Hermitian form, and 
$\Omega$ a 2-form on $M$ constructed
from $g$ as follows:
\begin{equation}\label{_Omega_Equation_}
\Omega := g(J\cdot, \cdot) +\1 g(K\cdot, \cdot)
\end{equation}
Then, $\Omega$ is a $(2,0)$-form on $(M,I)$ 
as an elementary linear-algebraic argument implies
(\cite{_Besse:Einst_Manifo_}).

 The hyperk\"ahler
condition can be written down as $d\Omega=0$ 
(\cite{_Besse:Einst_Manifo_}).
The HKT condition is weaker:

\hfill

\definition\label{_HKT_Definition_}
A quaternionic Hermitian metric is called an HKT-metric if
\begin{equation}\label{_HKT_intro_Equation_}
\6(\Omega)=0,
\end{equation}
where $\6:\; \Lambda^{2,0}_I(M) \arrow \Lambda^{3,0}_I(M)$
is the Dolbeault differential on $(M, I)$,
and $\Omega$ the $(2,0)$-form on $(M, I)$ constructed
from $g$ as in \eqref{_Omega_Equation_}. 

\hfill

In most examples a hypercomplex manifold
admits an HKT metric, except a few cases. For a 
long time, it was conjectured that any compact hypercomplex manifold
admits an HKT metric. However, in \cite{_Fino_Gra_}  
A. Fino and G. Grantcharov found an example of a compact
hypercomplex manifold not admitting an HKT-metric.

\hfill

The (2,0)-form $\Omega$ defined in \eqref{_Omega_Equation_}
is nowhere degenerate for any quaternionic Hermitian metric
$g$. This gives a nowhere degenerate section
$\Omega^{n}\in K(M,I)$ of the canonical
bundle, $n=\dim_{\Bbb H}M$. This section is in 
general non-holomorphic.

Clearly, the real structure $\eta \arrow J(\bar\eta)$ 
of \eqref{_real_stru_from_J_Equation_} preserves
$\Omega$, and therefore, preserves
$\Theta:=\Omega^{n}\in K(M,I)$.
Consider the 1-forms
\[
\theta:= \frac{\nabla^{1,0}_K\Theta}{\Theta}, \ \ 
\bar\theta:= \frac{\nabla^{0,1}_K\Theta}{\Theta}.
\]
where $\nabla_K$ is the Obata connection on $K(M,I)$.
If one considers the trivialization of $K(M,I)$
defined using $\Theta$, then $\theta+\bar\theta$
is the connection 1-form of the Obata connection
$\nabla_K$ on $K(M,I)$:
\[
\nabla_K = \nabla_{triv} + \theta +\bar\theta
\]
where $\nabla_{triv}$ is a trivial connection on
$K(M,I)$ preserving $\Theta$.
By writing also
\[
\nabla_{K^{1/2}} = \nabla_{triv} + \frac{1}{2}\theta +\frac{1}{2}\bar\theta,
\]
we obtain a connection on a square root $K^{1/2}$ of $K(M,I)$.
This way, the holomorphic line bundle $K^{1/2}$ is equipped with a connection
$\nabla_{K^{1/2}}$ and a real structure.

We equip $K^{1/2}$ also with a Hermitian structure,
in such  a way that the standard section
$\sqrt \Theta\in K^{1/2}$ has unit norm. The connection
$\nabla_{K^{1/2}}$ (which we also call the Obata connection)
is, generally speaking, non-Hermitian. 

Let $\bar\6:= \nabla^{0,1}_{K^{1/2}}$ be 
the $(0,1)$-part of $\nabla_{K^{1/2}}$, and
$\bar\6_J:= -J \circ \nabla^{1,0}_{K^{1/2}} \circ J$
be $\nabla^{1,0}_{K^{1/2}}$ twisted by $J$. 
The operator $J$ anticommutes with $I$, and, therefore,
maps $(p,q)$-forms on $(M,I)$ to $(q,p)$-forms. 
Therefore, both $\bar\6_J$ and $\bar\6$ 
map $\Lambda^{p,q}_I(M)\otimes K^{1/2}$ to
$\Lambda^{p,q+1}_I(M)\otimes K^{1/2}$. It is easy to check
that 
\[
\Lambda^{p,q}_I(M)\otimes K^{1/2}\stackrel{\bar\6, \bar\6_J}\arrow
\Lambda^{p,q+1}_I(M)\otimes K^{1/2}
\]
is a bicomplex
(see \cite{_Verbitsky:HKT_}), that is, $\bar\6$ and $\bar\6_J$
anticommute and square to zero.

In \cite{_Verbitsky:HKT_},
the following theorem was proven. 

\hfill

\theorem\label{_HKT_SS_Theorem_}
Let $M$ be a compact HKT-manifold, and
\[
\Lambda^{0,p}_I(M)\otimes K^{1/2}\stackrel{\bar\6, \bar\6_J}\arrow
\Lambda^{0,p+1}_I(M)\otimes K^{1/2}
\]
the bicomplex of $K^{1/2}$-valued $(0,p)$-forms
constructed above. Denote by $\Delta_{\bar\6}$, 
$\Delta_{\bar\6_J}$ the Laplacians on $\Lambda_I^{0,p}(M)\otimes K^{1/2}$
associated with the Hermitian structure defined above
$\Delta_{\bar\6}:= \bar\6 \bar\6^* + \bar\6^* \bar\6^*$,
$\Delta_{\bar\6_J}:= \bar\6_J \bar\6_J^* + \bar\6_J^* \bar\6_J^*$.
Then 
\begin{description}
\item[(i)]  The holomorphic cohomology of $K^{1/2}$ is naturally
identified with $\Delta_{\bar\6}$-harmonic forms
\item[(ii)]  $\Delta_{\bar\6}=\Delta_{\bar\6_J}$
\item[(iii)] Consider the operator 
\[ L_{\bar\Omega}:= \Lambda^{0,p}_I(M)\otimes K^{1/2}\arrow 
   \Lambda^{0,p+2}_I(M)\otimes K^{1/2},
\]
$\eta \arrow \eta \wedge \bar\Omega$, where
$\bar\Omega$ is complex conjugate to the HKT-form $\Omega$.
Denote by $\Lambda_{\bar\Omega}$ the Hermitian adjoint
operator, and let 
\[ H:\; \Lambda^{0,p}_I(M)\otimes K^{1/2}\arrow 
   \Lambda^{0,p}_I(M)\otimes K^{1/2}
\]
be a scalar operator acting on $(0,p)$-forms
as a multiplication by $(n-p)$, $n = \dim_{\Bbb H}M$.
Then $\langle L_{\bar\Omega}, \Lambda_{\bar\Omega}, H\rangle$
form an $\goth{sl}(2)$-triple in 
$\End\left(\Lambda^{0,*}_I(M)\otimes K^{1/2}\right).$
Moreover, these operators commute with 
$\Delta_{\bar\6}$.
\end{description}

{\bf Proof:} See \cite{_Verbitsky:HKT_}. \endproof

\hfill

\ref{_HKT_SS_Theorem_} gives a natural $\goth{sl}(2)$-action
on the holomorphic sheaf cohomology of $K^{1/2}$. This result
is analogous to the Lefschetz theorem about $\goth{sl}(2)$-action
on the cohomology of a compact K\"ahler manifold, and it is
proven in a similar way.

\hfill

We use \ref{_HKT_SS_Theorem_} to obtain the following result. 

\hfill

\theorem \label{_K_triv=>SL(n,H)_Theorem_}
Let $(M,I,J,K)$ be a compact hypercomplex manifold admitting
an HKT-metric, $\dim_{\Bbb H}M=n$, and $K(M,I)= \Lambda^{2n,0}_I(M)$
the canonical bundle of $(M,I)$. Assume that $K(M,I)$
is trivial as a holomorphic line bundle. Then the
Obata holonomy of $M$ is contained in $SL(n,{\Bbb H})$.

\hfill

{\bf Proof:} Let $\eta$ be a nowhere vanishing holomorphic section
of $K(M,I)$. As \ref{_SL_holo_Claim_} implies,
to prove \ref{_K_triv=>SL(n,H)_Theorem_},
we need to show that 
\begin{equation}\label{_holo_secti_para_Equation_}
\nabla_K(\eta)=0,
\end{equation}
where $\nabla_K$ is the connection on $K(M,I)$
induced by the Obata connection. Let $\eta_0:= \sqrt\eta$
be the corresponding holomorphic section of $K^{1/2}= \sqrt{K(M,I)}$.
Then 
\[ 
\frac{\nabla_K(\eta)}{\eta} = 2 \frac{\nabla_{K^{1/2}}(\eta_0)}{\eta_0},
\] 
hence \eqref{_holo_secti_para_Equation_} would follow from
\begin{equation}\label{_holo_secti_para_for_root_Equation_}
\nabla_{K^{1/2}}(\eta_0)=0. 
\end{equation}
Since $\eta_0$ is holomorphic, $\Delta_{\bar\6}(\eta_0)=0$,
and therefore $\Delta_{\bar\6_J}(\eta_0)=0$ and $\bar\6_J(\eta_0)=0$.
Since \[ \bar\6_J= -J \circ \nabla^{1,0}_{K^{1/2}} \circ J, \]
$\bar\6_J(\eta_0)=0$ implies $\nabla^{1,0}_{K^{1/2}}(\eta_0)=0$.
We obtain that
\[
\nabla_{K^{1/2}}(\eta_0) = \nabla^{1,0}_{K^{1/2}}(\eta_0)+ \bar\6(\eta_0)=0
\]
This proves \eqref{_holo_secti_para_for_root_Equation_}
and \eqref{_holo_secti_para_Equation_}. We proved 
\ref{_K_triv=>SL(n,H)_Theorem_}. \endproof


\section[Appendix. Hypercomplex 
manifolds with holonomy in $SL(n, {\Bbb H})$]{Appendix. Hypercomplex 
manifolds with \\ holonomy in $SL(n, {\Bbb H})$}


\subsection{Irreducible holonomy and hypercomplex geometry}

Let $(M, \nabla)$ be a manifold with torsion-free
connection in $TM$, and $\Hol_0(\nabla)$ its local holonomy
group. One is interested in irreducible
holonomies, that is, the groups $G=\Hol_0(\nabla)$ 
such that the action of $G$ on $TM$ is irreducible.
Classification of all Lie groups occuring this way
has a long history, starting from the works of Eli Cartan
in 1920-ies. In 1955 M. Berger published a paper 
\cite{_Berger:holonomies_}, which
contained a list of irreducible holonomies, for Levi-Civita
connections and for general torsion-free
affine connections. This
work opened a new chapter in the study
of differential geometry of Riemannian
manifolds. The Berger's classification 
of Riemannian holonomies is one of the
cornerstones of modern differential geometry,
strongly influencing physics and 
algebraic geometry, and this 
influence increases still.

Berger's list of non-Riemannian holonomies was
largely ignored. In fact, it took
almost 40 years until the omissions 
in Berger's list were found by 
S. Merkulov and L.  Schwachh\"ofer,
who provided a complete 
classification of irreducible holonomies
of torsion-free affine connections
in \cite{_Merkulov_Sch:long_}.

The hypercomplex manifolds are equipped
with the Obata connection, which is also tor\-sion-\-free.
The Obata connection preserves a quaternionic action, and therefore,
the holonomy of the Obata connection belongs to $GL(n,{\Bbb H})$.
Conversely, if $(M,\nabla)$ is equipped with a torsion-free
affine connection with (global) holonomy in $GL(n,{\Bbb H})$,
then $M$ is hypercomplex. 

For many (or most of) 
the groups from Merkulov-Schwachh\"ofer
list, it is not clear whether they can be realized
as holonomies of {\bf compact} manifolds. It is not
clear which subgroups of $GL(n,{\Bbb H})$
can occur as local holonomies
of compact manidolds.

We do not know any example 
of a compact manifold with local holonomy
in $SL(n,{\Bbb H})$. In fact, all known
compact manifolds with local holonomy
in $SL(n,{\Bbb H})$ are locally hyperk\"ahler.

In this Appendix, we show that 
there exists a compact, hypercomplex, 
non-\-hyp\-er\-k\"ahler 
manifold with (global) holonomy in $SL(n, {\Bbb H})$.
We use the examples of hypercomplex structures on
nilmanifolds constructed 
in \cite{_Barberis_Dotti_} by M. L. Barberis
and I. Dotti, and much studied since then
by M. L. Barberis, I. Dotti and A. Fino.

\subsection{Abelian hypercomplex structures on nilmanifolds}

\definition\label{_abelian_hc_Definition_}
Let $\g$ be a Lie algebra over $\R$, and
${\Bbb H}\otimes \g \arrow \g$ be an
action of the algebra of quaternions on $\g$. 
This action is called {\bf abelian} if 
\[
[Ix,Iy] = [Jx, Jy] = [Kx,Ky] = [x,y]
\]
for all $x, y \in \g$. 
Abelian quaternionic structures on Lie algebras 
were studied in \cite{_Barberis_Dotti_} and
\cite{_Dotti_Fino:8-dim_}. 

If $\g$ admits such an action, then
$\g$ is solvable (\cite{_Dotti_Fino:Bilbao_}).

In \cite{_Dotti_Fino:8-dim_} the 8-dimensional
Lie algebras with an abelian quaternionic 
structure were classified. 

\hfill

In a similar way one defines an abelian complex
structrure on a Lie algebra: the complex structure
$I:\; \g\arrow \g$ is abelian if $[Ix,Iy] = [x,y]$
for all $x, y \in \g$. The following lemma is elementary

\hfill

\lemma\label{_abelian_comple_on_Lie_alg_Lemma_}
Let $\g$ be a Lie algebra over $\R$ equipped with
an abelian complex structure $I$. Denote by $G$ the corresponding
Lie group, and let $\cal I$ be the left-invariant almost complex
structure operator on $G$ obtained by transporting $I$ around.
Then $\cal I$ is integrable.

\hfill

{\bf Proof:} To show that ${\cal I}$ is integrable, we need to check
that 
\begin{equation}\label{_integrabi_Equation_}
[x,y] \in T^{1,0}_{\cal I}G \text{ \ for any \ } x, y \in T^{1,0}_{\cal I} G.
\end{equation}
This condition is $C^\infty$-linear, hence we may assume that
$x$, $y$ are left-\-in\-va\-ri\-ant. The commutator of left-\-in\-va\-ri\-ant
vector fields is determined by the Lie algebra. Therefore,
\eqref{_integrabi_Equation_} is implied by the same equation
in the Lie algebra:
\[
[x,y] \in \g^{1,0}_I \text{ \ for any \ } 
x, y \in \g^{1,0}_I\subset \g\otimes \C
\]
We write $x= u+\1 Iu$, $y=v+\1 Iv$, where 
$u,v \in \g$. Then
\begin{align*} 
  [x, y] & = [u, v] - [Iu, Iv] + \1 [u, Iv] + \1 [Iu, v] \\ & 
   = [u, v] - [Iu, Iv] + \1 [u, Iv] - \1 [u, Iv] =0
\end{align*}
(the middle equality is true because $I$ is abelian).
This proves \ref{_abelian_comple_on_Lie_alg_Lemma_}.
\endproof

\hfill

The same argument also implies the following claim

\hfill

\claim\label{_commu_1,0_zero_Claim_}
Let $\g$ be a Lie algebra with an abelian complex structure
and $x,y\in \g^{1,0}\subset \g\otimes \C$
the $(1,0)$-vectors in its complexification.
Then $[x,y]=0$.
\endproof

\hfill

\definition
Let $g$ be a nilpotent  Lie algebra equipped
with an abelian quaternionic algebra action. 
Using \ref{_abelian_comple_on_Lie_alg_Lemma_},
we can extend this action to a left-invariant
hypercomplex structure on the corresponding
Lie group $G$. For any discrete subgroup $\Gamma\subset G$,
$M:= \Gamma\backslash G$ is also hypercomplex. We call
$M$ {\bf a nilmanifold equipped with an
abelian hypercomplex structure}. Many 
compact examples of such nilmanifolds are 
known, see e.g. \cite{_Dotti_Fino:8-dim_}.

\hfill

In \cite{_Dotti_Fino:HKT_} (see also
\cite{_Fino_Gra_}), an HKT-metric was constructed
on any nilmanifold equipped with an
abelian hypercomplex structure. We 
give a short version of this construction,
using the description of HKT-metrics
given in \cite{_Verbitsky_HKT_exa_} 
and \cite{_Verbitsky:HC_Kahler_}.

\hfill

Let $(M, I, J, K)$ be a hypercomplex manifold. Since $J$ and $I$
anticommute, $J$ maps $(p,q)$-forms on $(M,I)$
to $(q,p)$-forms:
\[ J:\; \Lambda^{p,q}_I(M) \arrow \Lambda^{q,p}_I(M).
\]

\definition
Let $\eta\in \Lambda^{2,0}_I(M)$ be a (2,0)-form on $(M,I)$.
Then $\eta$ is called {\bf $J$-real} if $J(\eta) = \bar\eta$,
and {\bf $J$-positive} if for any $x\in T^{1,0}(M,I)$,
$\eta(x, J(\bar x)) \geq 0$. We say that $\eta$ is 
{\bf strictly $J$-positive} if this inequality is strict
for all $x\neq 0$.

Denote the space of $J$-real, 
strictly $J$-positive  $(2,0)$-forms
by $\Lambda^{2,0}_{>0}(M,I)$.

\hfill

We need the following linear-algebraic lemma, which is
well known (its proof can be found e.g. \cite{_Verbitsky_HKT_exa_}).

\hfill

\lemma\label{_2,0_forms_and_q_H_metrics_Lemma_}
Let $M$ be a hypercomplex manifold. Then 
$\Lambda^{2,0}_{>0}(M,I)$ is in one-to-one correspondence
with the set of quaternionic Hermitian metrics $g$ on $M$.
This correspondence is given by
\[ g \arrow  g(J \cdot, \cdot)+ \1 g(K \cdot, \cdot),\]
and the inverse correspondence by
\begin{equation}\label{_g_from_Omega_Equation_}
\Omega \arrow g(x,y):= \Omega(x, J(\bar y)).
\end{equation}

\endproof

\hfill

From \ref{_2,0_forms_and_q_H_metrics_Lemma_}
and \ref{_HKT_Definition_},
it follows that to define an HKT-metric on $M$
it is  sufficient to find a $(2,0)$-form
$\Omega\in \Lambda^{2,0}_{>0}(M,I)$
satisfying $\6 \Omega=0$.

\hfill

Return now to the case of nilmanifolds. 
The vector space $\g$ is isomorphic to ${\Bbb H}^n$
as a quaternionic vector space, therefore,
it is possible to find a form $\eta\in \Lambda^{2,0}_I\g$
which is $J$-real and strictly $J$-positive. 
Extending $\eta$ to a left-invariant
form $M = \Gamma\backslash G$, we obtain
a $(2,0)$-form $\Omega\in \Lambda^{2,0}_{>0}(M,I)$.
Any left-invariant $(p,0)$-form $\eta$ satisfies
$\6\eta=0$ as follows from \ref{_commu_1,0_zero_Claim_}
and the Cartan's formula. This gives us that $\6 \Omega=0$.
We just proved the following claim.

\hfill

\claim\label{_nilma_HKT_Claim_}
\cite{_Dotti_Fino:HKT_}
Let $M$ be a nilmanifold equipped with
an abelian left-invariant hypercomplex structure.
Then $M$ admits an HKT-metric.

\endproof 

\hfill

\subsection{Canonical bundle of a nilmanifold}

\proposition\label{_can_class_trivi_Proposition_}
Let $\g$, $\dim_\R\g =2n$ 
be a nilpotent Lie algebra over $\R$ equipped with an
abelian complex structure, and $G$ the corresponding
Lie group. Using \ref{_abelian_comple_on_Lie_alg_Lemma_},
we may consider $G$ as a complex manifold.
Let $\Theta$ be a left-invariant section
of the canonical bundle $\Theta \in \Lambda^{n,0}(G)$.
Then 
\begin{description}
\item[(i)] $\Theta$ is closed: $d\Theta=0$
\item[(ii)] $\Theta$ is holomorphic. 
\end{description}
\hfill

{\bf Proof:} Let $g_1, ..., g_{2n}$ be a basis in $\g$,
chosen in such a way that 
\begin{equation}\label{_c^k_ij_Equation_}
[g_i, g_j] = \sum c^k_{i,j}g_k, \ \ 
c^k_{i,j}=0 \text{\ for \ } k\leq \max(i,j).
\end{equation}
(such a basis always exists because $\g$ is nilpotent).
Let $g_i^*$ be the dual basis, and $\xi_i\in \Lambda^1 G$
be the corresponding basis in the space of left-invariant differential
forms. Clearly, Cartan's formula implies
\begin{equation}\label{_d_xi_i_Cartan_Equation_}
d \xi_k = \sum c^k_{i,j} \xi_i \wedge \xi_j,  \ \ 
  c^k_{i,j}=0 \text{\ for \ } k\leq \max(i,j).
\end{equation}
Since the complex structure in $\g$ is abelian,
we may always chose a basis $g_1, ..., g_{2n}$ 
such that $I(g_{2i-1})= g_{2i}$ and 
\eqref{_c^k_ij_Equation_} is satisfied. 
Let $h_i:= g_{2i-1}+ \1 g_{2i}$
be the corresponding basis in $\g^{1,0}$, and
$\theta_i\in \Lambda^{1,0} G$ the dual basis
of left-invariant $(1,0)$-forms on $G$.
Then $\bigwedge \theta_i$ is a non-trivial
section of the canonical bundle
$\Lambda^{n,0}(G)$. Being left-invariant,
this section is
unique up to a constant:
\[ \Theta = c \bigwedge \theta_i.\]

\hfill

Write 
\[
d \theta_k = \sum m^k_{i,j} \theta_i \wedge \bar\theta_j,    
\]
Equation \eqref{_d_xi_i_Cartan_Equation_} implies that
$m^k_{i,j}=0$ for  $k\leq \max(i,j)$.
Then 
\begin{equation}\label{_d_of_section_cano_Equation_}
\begin{aligned}
d \left(\bigwedge \theta_i\right) = &
\sum_k \bigg [\theta_1\wedge \theta_2 \wedge ...
\wedge\theta_{k-1} \wedge \theta_{k+1} \wedge ... \wedge \theta_n\\
&\wedge \sum m^k_{i,j} \theta_i \wedge \bar\theta_j\bigg]\\
m^k_{i,j} & =0 \text{\ for \ } k\leq \max(i,j)
\end{aligned}
\end{equation}
Every term of \eqref{_d_of_section_cano_Equation_}
in brackets clearly vanishes, because it necessarily involves
a product of $\theta_i$ and $\theta_i$, $i<k$, as
$m^k_{i,j}=0$ for  $k\leq \max(i,j)$.
We have thus proved \ref{_can_class_trivi_Proposition_} (i).
The second part of \ref{_can_class_trivi_Proposition_}
is a formal consequence of the first. 
Indeed, $\bar\6 \Theta = d\Theta$, because
$\6\Theta=0$ (there are no non-zero $(n+1,0)$-forms,
and $d = \6 + \bar\6$). This proves
\ref{_can_class_trivi_Proposition_}.
\endproof

\hfill

Comparing \ref{_nilma_HKT_Claim_}, 
\ref{_can_class_trivi_Proposition_} and 
\ref{_K_triv=>SL(n,H)_Theorem_}, 
we obtain the following corollary.

\hfill

\corollary\label{_SL(N,H)_nilma_Corollary_}
Let $M$ be a compact nilmanifold equipped with
an abelian left-invariant hypercomplex structure.
Then the holonomy of the Obata connection on $M$ 
is contained in $SL(n, {\Bbb H})$.

\endproof

\hfill

\remark
In the examples of nilmanifolds considered in
\cite{_Dotti_Fino:8-dim_}, the local holonomy
of Obata connection is trivial. Therefore,
the global holonomy is equal to the 
action of the fundamental group $\Gamma$,
induced by the monodromy of this flat connection.
This can be used to obtain 
\ref{_SL(N,H)_nilma_Corollary_} directly.

\hfill

{\bf Acknowledgements:}
This paper appeared as a result of a very rewarding
colloboration with S. Alesker. I am also grateful
to D. Kaledin and D. Kazhdan for interesting discussions.

\hfill

{\small

\hfill

\noindent {\sc Misha Verbitsky\\
University of Glasgow, Department of Mathematics, \\
15 University Gardens, Glasgow G12 8QW, Scotland}, \\
{\sc  Institute of Theoretical and
Experimental Physics \\
B. Cheremushkinskaya, 25, Moscow, 117259, Russia }\\
\tt verbit@maths.gla.ac.uk, \ \  verbit@mccme.ru 
}

\end{document}